\documentclass[12pt]{article}
\usepackage{amsfonts,hyperref,textcomp}
\usepackage{caption,subcaption}
\usepackage{graphicx}
\usepackage{blkarray}
\usepackage{array}
\usepackage[section]{placeins}
\usepackage{amssymb}
\usepackage{booktabs}
\usepackage{multirow}
\usepackage{float}
\usepackage{epstopdf}
\usepackage{color}
\usepackage{amsmath}
\usepackage[normalem]{ulem}
\usepackage{authblk}

\usepackage{enumitem}
\newlist{worddefs}{description}{1}
\setlist[worddefs]{font=\sffamily\bfseries, labelindent=\parindent, leftmargin=6em, style=sameline}

\newtheorem{theorem}{\bf Theorem}[section]

\begin{document}

%%%% Article title to be placed here
%% \title{Dynamics in simple delay differential equations with periodic mixed feedback}
\title{Periodic patterns in simple biological differential delay models}

\author[1,*]{A. Ivanov}
\author[2]{S. Shelyag}

%%%%%%%%% Insert author address here
\affil[1]{Department of Mathematics, Pennsylvania State University, Dallas, PA 18612, USA}
\affil[2]{College of Science and Engineering, Flinders University, Tonsley, Adelaide, SA 5042, Australia}
\affil[*]{Corresponding author A. Ivanov, \it{afi1@psu.edu}}

%%%% Keyword entries to be placed here %%%%
%\keywords{Delay differential equations, periodic coefficients, mixed positive-negative feedback, periodic solutions, reduction to interval maps, stability}

%%%% Insert corresponding author and its email address}
%\corres{S. Shelyag\\
%\email{sergiy.shelyag@flinders.edu.au}}

%\corres{A. Ivanov\\
%\email{afi1@psu.edu}}

\maketitle

%%%% Abstract text to be placed here %%%%%%%%%%%%
\begin{abstract}
Periodic patterns in dynamical behaviours of biological models described by simple form differential delay equations are studied. 
Mathematical models are given by a class of scalar delay differential equations with a multiplicative time periodic mixed coefficient and a nonlinear delayed negative 
feedback. The dynamics is studied analytically with supportive numerical simulation and justification of the theoretical outcomes.
The principal  nonlinearity involving the state variable is of the negative feedback type, the periodic multiplicative coefficient can change its sign, leading to equations with mixed positive-negative feedback. The existence of slowly oscillating periodic solutions of two different types is established. The theoretical analysis and derivation are based on the reduction of dynamics in the delay equations to that of interval maps. The theoretical outcomes are verified and supported by comprehensive numerical computations. The differential delay equations considered are generalisations of some well-known autonomous models from biological applications. 
\end{abstract}

\section{Introduction}

Differential delay equations are widely used as mathematical models of various real-world phenomena where the after-effects are the 
intrinsic features of their 
functioning. Such applications are broad and numerous; some distinctive samples can be found in fundamental monographs on the theory and 
applications of these equations. For comprehensive description of various models from applications, together with the relevant theoretical
expositions, one can consult monographs \cite{Ern09,GlaMac88,KolMys99,Kua93,Smi11}. For the basics and fundamentals of the theory the 
reader is referred to classical  monographs \cite{DieSvGSVLWal95,HalSVL93}.

Simple form nonlinear delay differential equations of the type 
\begin{equation}\label{DDE}
    \frac{dx(t)}{dt}=-\mu x(t)+f(x(t-\tau)), \quad \mu>0, \tau>0
\end{equation}
serve as mathematical models in a variety of applications, in particular in biological sciences \cite{GlaMac88,Had79,Kua93,Smi11}. 
Of interest are population models which include the well-known Nicholson's blowfly equation 
\cite{GurBlyNis80,PerMalCou78} and Mackey-Glass physiological models \cite{MacGla77,GlaMac88}.
Such models exhibit diverse behaviours of solutions including the global asymptotic stability of equilibria, existence of periodic solutions with 
different types of stability, transition to complex dynamics, and chaotic behaviours. The theoretically derived dynamics in many cases closely correspond 
to observed  phenomena in corresponding biological models.

Under the negative feedback assumption in equation (\ref{DDE}) and the instability of the unique positive equilibrium, the periodic solutions typically 
exist. This a phenomenon cannot be adequately described by corresponding ordinary differential equations, where only stable or unstable equilibria may 
exist. Therefore, the presence of delays in the differential models plays determining role in the existence of periodicity and more complex dynamics
(e.g., chaos). The complex dynamics in models of type (\ref{DDE}) is observed when it has the typical properties of physiological systems, - 
a combination of the positive and negative feedback with respect to equilibria \cite{GlaMac88,MacGla77}. In the case of negative equilibrium the state's deviation from
the  equilibrium in one direction forces its movement in the opposite direction. In the case of the positive feedback such a deviation is enhanced and its amplitude increases. 

In view of the temporal periodicity factors in the environment, such as seasonal changes and circadian rhythms among others,  
a more adequate delay differential equation to describe the processes has the form 
\begin{equation}\label{DDE-per}
    \frac{dx(t)}{dt}=-\mu(t) x(t)+a(t) f(x(t-\tau(t))), 
\end{equation}
where the functions $\mu(t)>0, \tau(t)>0$ and $a(t)$ are all $T$-periodic.

The negative feedback condition in the models (\ref{DDE}) and (\ref{DDE-per}) is expressed by the inequality $x\cdot f(x)<0, x\ne0$, where $x(t)\equiv0$
is the (unique) equilibrium. 
The initially positive equilibrium $x_*>0$ in models (\ref{DDE})-(\ref{DDE-per}) can be translated to the zero equilibrium by a shift in the variable $x$. 
The positive feedback about the equilibrium $x(t)\equiv0$ is analytically described by $x\cdot f(x)>0, x\ne0.$ In this work we are interested in the 
dynamics of equation (\ref{DDE-per}) when both the positive and negative feedback can be present about the same equilibrium. Such a phenomenon  can be  achieved in 
DDE (\ref{DDE-per}) when $f(x)$ has one fixed type of the feedback, however, the multiplicative coefficient $a(t)$ changes sign. Under changing assumption, we show the existence of stable periodic solutions to equation (\ref{DDE-per}).

A well known and classical by now approach to discover particular dynamics in DDEs (\ref{DDE}) and (\ref{DDE-per}) is the use of piecewise constant functions 
$f, \mu$ and $a$. We refer the reader to the pioneering paper \cite{Pet80}, the follow up analyses on comlpicated (chaotic) behaviours in papers 
\cite{Wal81b}, and a review paper  \cite{IvaSha92}. In this work, we use a similar approach through the piecewise constant functions $\mu, a, f$ in model (\ref{DDE-per}). We start with a simplified version of the equation when $\mu(t)\equiv0.$

Consider first the delay differential equation
\begin{equation}
\label{eq1}
\frac{dx}{dt} = a(t)\,f(x(t-\tau)), \qquad f(x) = -\operatorname{sign}(x),
\end{equation}
where $\tau > 0$ is a fixed constant delay, $a(t)$ is a piecewise-constant function that alternates between positive and negative values, and $f(x)$ represents 
a  regulatory response with negative feedback. By explicit calculation we prove the existence of stable slowly oscillating periodic solutions of 
two types with period $T$ and $2T$, respectively. This is done in Section \ref{Results}, parts (b) and (c). In part (d) of the section we show that the stable periodic 
solutions persist when we make the piecewise constant functions continuous, by "smoothing" them in a small neighbourhood of the jump discontinuity points. 
We also provide some basic related preliminaries on delay differential equations (Section \ref{Basics}) as well as discuss possible further directions of 
research (Section \ref{Concl}). This work is based on and uses some of our previous approaches and results obtained in papers 
\cite{IvaLWShel2024,IvaSha92,IvaShel2023,IvaShel2024,IvaShel2024_MCA,IvaShel2025}. It further extends them to the case of delay differential models with the feedback of mixed type 
(negative-to-positive).

%%
%% Section "Methods"
%%
\section{Methods}\label{Methods}

In this work we use a classical by now approach of studying dynamical properties of delay differential equations, which can formally be viewed as infinite-dimensional dynamical systems, by reducing some of their dynamics to finite-dimensional discrete maps. In the simplest case of scalar equations of the form 
(\ref{DDE}) and (\ref{DDE-per}) such a reduction can be done by approximating the nonlinearity $f$ by a piece-wise constant one, and reduction of the dynamics on a subset of its 
phase space to one-dimensional (interval) maps. The pioneering work here belongs to Peters \cite{Pet80} and Walther \cite{Wal81b} together with others whose 
papers follow shortly thereafter \cite{AdHM,AdHW,HalLin86}. The paper \cite{IvaSha92} is a review work on a particular equation of the form (\ref{DDE}) together with a list of additional references on the subject. 

We follow the basic ideas of the mentioned above papers by applying them to the simplest form DDE (\ref{DDE-per}). We focus on the so-called slowly oscillating solutions which are characterised by the fact that the distance between consecutive zeros is greater than the delay. Though the negative feedback property of the function $f$ and the variable sign of the multiplicative coefficient $a$ make the feedback in the system of the alternative positive-negative nature, we manage to identify the parameters' range for a subset of initial functions which result in forward slowly oscillating solutions. This reduces the dynamics on the selected subset of the phase space $C$ to the dynamics of an interval map on the corresponding initial set of functions, which are now characterised by a single parameter. The fixed points of the interval map correspond to the periodic solutions of the differential delay equation (\ref{DDE-per}). The attracting or repelling nature of the fixed points translates to the same type of stability of the periodic solutions.

\section{Basic assumptions and preliminaries}\label{Basics}
The standard phase space for DDE (\ref{DDE-per}) with continuous functions $\mu, f, a$ is the set $C([-\tau_0,0], \mathbb{R}):=C$ of continuous real-valued functions
on the initial interval $[-\tau_0,0]$ , where $\tau_0=\max\{\tau(t), t\in[0,T]\}$. With an initial function $\varphi\in C$ the corresponding solution 
$x(t)=x_{\varphi}(t)$ to equation (\ref{DDE-per}) is unique and exists for all $t\ge 0$ \cite{HalSVL93,DieSvGSVLWal95}.

The same set of initial functions can also be chosen to construct the forward solutions of discontinuous DDEs, in particular those considered in Section 
\ref{Results}. Since both functions $f$ and $a$ can only have jump discontinuities, under a solution $x_{\varphi}(t), t\ge0,$ of the corresponding 
equation (\ref{DDE-per}) we mean a continuous function satisfying it for all $t\ge0$, except a discrete set of isolated $t$-values (see, e.g., Figure~\ref{fig1} and~\ref{fig2} 
in Section \ref{Results}). At those points the one-sided derivatives of the solution exist but are of different values, making a solution of the ``corner-type'' shape.
In order to have such a solution defined for $t\ge0$ one must impose additional restrictions on the set of initial functions. Namely, the initial functions 
must intersect the discontinuity values for $f$ and $a$ transversely. In our considerations in Section \ref{Results} the initial functions are chosen such that the 
corresponding solutions exist for all forward time.

%\textcolor{red}{... describe appropriately basics related to \textit{discontinuous} initial data/functions ...}

%%
%% Section 4 "Results"
%%
\section{Results}\label{Results}

In this section we show the existence of slowly oscillating periodic solutions to DDE (\ref{DDE-per}) with the same or the double of the period $T$ as the 
period of the coefficient $a(t)$. We start  with the case when both the nonlinearity $f$ and the coefficient $a$ are piecewise constant. Those calculations
are more involved but analogous to our previous results derived in papers \cite{IvaShel2024}. We then apply a standard smoothing procedure to make the 
functions involved continuous. In fact, one can achieve any degree of smoothness of $a$ and $f$ by making appropriate matching connections in a small 
neighbourhood of the discontinuity points. We only describe the procedure in sufficient detail of our case; we shall use the existing results of ours 
derived in \cite{IvaShel2024_MCA,IvaLWShel2024} for a more general equation of the type (\ref{DDE-per}).

%%
%% subsection (a)
%%
\subsection{Piece-wise constant nonlinearities}\label{PCN}
In this subsection we consider the particular case of equation (\ref{DDE-per}) with the decay coefficient $\mu(t)\equiv0$ and the constant delay $\tau=1$: 
\begin{equation}\label{P-C dde}
  x^\prime(t)=a(t)f(x(t-1)),\qquad t\ge0.
\end{equation}
The multiplicative coefficient $a(t)$ is $T$-periodic, with the period greater than the delay, $T>1$. The nonlinearity $f$ is of the negative 
feedback type, $x f(x)<0, x\ne0.$ The coefficient $a(t)$ is allowed to be sign changing on the periodic interval $[0,T]$.
Note that the case of general delay $\tau>0$ can always be normalised to $\tau=1$ by time rescaling $t=\tau\cdot s$.

We start with the case when the nonlinearity $f$ is the negative sign function    
$$
f(x)=f_0(x)=-\text{sign} (x)=\begin{cases}
+1\;\text{if}\; x<0\\
{\;\;\;}0\;\text{if}\; x=0\\
-1\;\text{if}\; x>0,
\end{cases}
$$ 
and the $T$-periodic coefficient $a(t)$ is a piecewise constant function defined by three constants $a_1>0, a_2>0$ and $a_3>0$ on respective time 
segments which lengths are given by constants $p_1>0,p_2>0,$ and $ p_3>0$  as follows
$$
a(t)=A_0(t)=\begin{cases}
a_1\;\text{if}\; t\in[0,p_1)\\
a_2\;\text{if}\; t\in[p_1,p_1+p_2)\\
-a_3\;\text{if}\; t\in[p_1+p_2,p_1+p_2+p_3)\\
\text{periodic extension on}\; \mathbb R\;\text{outside}\\ \text{interval}\; [0,T), T=p_1+p_2+p_3,
\end{cases}
$$
where $a_1, a_2, a_3, p_1, p_2, p_3$ are all positive constants.

%%
%% Type I PSs
%%
\subsection{Periodic solutions of type I}\label{PS-I}

Given $f=f_0$ and $a=A_0$ as above, for arbitrary initial function $\phi\in C$ the corresponding solution $x=x(t,\phi)$ is
constructed for all forward times $t\ge0$. It is a continuous piecewise differentiable function composed of consecutive matched segments 
of affine (linear) functions.

\begin{figure*}
    \centering   \includegraphics[width=0.9\textwidth]{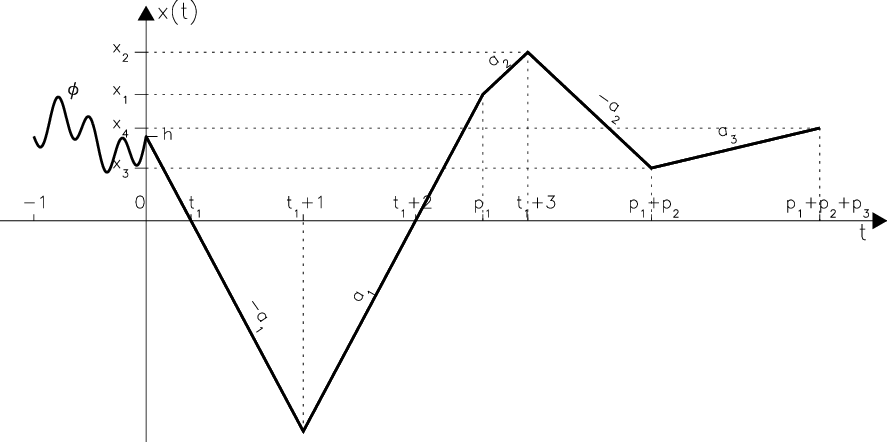}
    \caption{Slowly oscillating piecewise affine solution (Type I) on initial periodic interval $[0,T], T=p_1+p_2+p_3.$}
    \label{fig1}
\end{figure*}

We assume that periodic solutions have the form as depicted in Figure~\ref{fig1} (this will be shown to be the case for particular values of the
parameters $a_1,a_2,a_3,p_1,p_2,p_3$). For arbitrary initial function $\phi\in C$ such that $\phi(s)>0\; \forall s\in[-1,0],$ 
the corresponding
forward solution $x(t),t\ge0,$ depends only on the value $\phi(0):=h>0$ and does not depend on the other values $\phi(s), s\in[-1,0)$
of the initial interval. There exist consecutive zeros $t_1>0$ and $t_2=t_1+2$ of the solution such that $x(t)>0, t\in[0,t_1)$, 
$x(t)<0, t\in(t_1,t_1+2)$, and $x(t)>0, t\in[t_1+2,p_1]$. The necessary condition for such first two zeros to exist is $p_1>2$.
We also assume that $p_1<t_1+3<p_1+p_2$ and that the solution remains positive on the interval $[p_1,p_1+p_2]$. For the latter to be 
true one only has to require that $x_3:=x(p_1+p_2)>0$. Since the solution's slope on the interval $[p_1+p_2,p_1+p_2+p_3]$ is positive 
$a_3>0$ the assumption $x_3>0$ implies that the entire segment of the solution $x(t)$ is positive on the interval 
$(t_1+2, p_1+p_2+p_3)$. 

Inspired by our previous results in \cite{IvaShel2024_MCA} and numerical simulations of our particular DDE (\ref{eq1}), we are looking for periodic
solutions of the form shown in Figure \ref{fig1}.
%%
%% Fig.1
%%

We next calculate the explicit values of $t_1, x_1=x(p_1), x_2=x(t_1+3), x_3=x(p_1+p_2), x_4=x(p_1+p_2+p_3)$ in terms of the parameters $a_1,a_2,a_3, p_1, p_2, p_3$ 
and given initial value $h>0$. It is straightforward to find that:
$$
t_1=\frac{h}{a_1},\quad x_1=x(p_1)=-h+a_1p_1-2a_1,\; 
$$
$$
x_2=x(t_1+3)=\left(\frac{a_2}{a_1}-1\right)h+a_1p_1-2a_1+3a_2-a_2p_1,
$$
\begin{align*}
    x_3=&x(p_1+p_2)=\left(2\frac{a_2}{a_1}-1\right)h+a_1(p_1-2)+\\
        &a_2[6-(2p_1+p_2)]
\end{align*}

\begin{align*}
    x_4=& x(p_1+p_2+p_3)=x_3+a_3p_3=\left(2\frac{a_2}{a_1}-1\right)h+\\
        &a_1(p_1-2)+a_2[6-(2p_1+p_2)]+a_3p_3:=mh+b.
\end{align*}

In addition to the previous assumptions the condition $b>0$, that is, 
\begin{equation}\label{b}
    b=a_1(p_1-2)+a_2[6-(2p_1+p_2)]+a_3p_3>0,
\end{equation}
guarantees that $x_3>0, x_4>0$, and thus the solution $x(t), t\in[0,p_1+p_2+p_3],$ is of desired shape (as shown in Figure~\ref{fig1}).

A fixed point $h_*>0$ of the affine map $F(x)=mx+b$ gives rise to a slowly oscillating periodic solution $x=x(t,h_*)$ of equation 
(\ref{DDE-per}). Moreover, by the linearity of $F$ and the construction, such periodic solution is asymptotically stable if $|m|<1$. 
The latter is equivalent to $0<a_2<a_1$.
The unique fixed point $x=h_*$ is easily found as 
\begin{equation}\label{h}
    h_*=\frac{b}{1-m}=\frac{a_1[a_1(p_1-2)+a_2[6-(2p_1+p_2)]+a_3p_3]}{2(a_1-a_2)}.
\end{equation}
Therefore, we arrive at the following
%%
%% Theorem 4.1
%%
\begin{theorem}\label{Thm1}
Suppose that the parameters $a_1,a_2,a_3,p_1,p_2,p_3$ are such that the inequality (\ref{b}) is satisfied and $a_1>a_2$. 
Then DDE (\ref{DDE-per}) has an asymptotically stable slowly oscillating periodic solution of period $T=p_1+p_2+p_3$. The periodic 
solution is generated by the  initial function $\phi(s)\equiv h_*, s\in[-1,0],$ where $h_*$ is given by (\ref{h}).
\end{theorem}

Note that due to the symmetry property of the nonlinearity $f$ (oddness) and the procedure of construction of the periodic solution
$x_{h_*}$, under the assumption of Theorem \ref{Thm1} there also exists the symmetric to $x_{h_*}$ periodic solution generated by
the initial function $\phi(s)\equiv -h_*$. The two periodic solutions are related by $x_{-h_*}(t)\equiv-x_{h_*}(t)$.

Some examples of the parameter values for which type-I solutions of (1.2) have been obtained numerically are given in Table \ref{table_single_T}.

\begin{table*}[h!]
\centering
\begin{tabular}{|c|c|c|c|c|c|c|c|}
\hline
$a_1$ & $a_2$ & $a_3$ & $p_1$ & $p_2$ & $p_3$ & $h^*$ & $T$ \\
\hline
0.5 & 0.1 & -0.1 & 3.0 & 1.0 & 0.5 & 0.2812 & 4.5 \\
0.5 & 0.1 & -0.5 & 3.0 & 1.0 & 0.5 & 0.4062 & 4.5 \\
1.0 & 0.1 & -0.1 & 3.0 & 1.0 & 1.0 & 0.5556 & 5 \\
5.0 & 0.1 & -1.0 & 3.0 & 1.0 & 1.0 & 3.01 & 5 \\
5.0 & 1.0 & -2.0 & 3.0 & 1.0 & 1.0 & 3.75 & 5 \\
5.0 & 3.0 & -1.0 & 3.0 & 1.0 & 1.0 & 3.75 & 5 \\
10.0 & 0.1 & -0.5 & 3.0 & 1.0 & 1.0 & 5.253 & 5 \\
10.0 & 0.5 & -0.5 & 3.0 & 1.0 & 1.0 & 5.263 & 5 \\
3.0 & 0.5 & -0.1 & 3.0 & 1.0 & 5.0 & 1.8 & 9 \\
5.0 & 3.0 & -0.1 & 3.0 & 1.0 & 5.0 & 3.125 & 9 \\
10.0 & 0.5 & -0.1 & 3.0 & 1.0 & 5.0 & 5.263 & 9 \\
10.0 & 3.0 & -0.1 & 3.0 & 1.0 & 5.0 & 5.357 & 9 \\
3.0 & 0.1 & -1.0 & 5.0 & 0.5 & 5.0 & 7.009 & 10.5 \\
0.5 & 0.1 & -2.0 & 5.0 & 1.0 & 0.5 & 1.25 & 6.5 \\
1.0 & 0.5 & -2.0 & 5.0 & 1.0 & 1.0 & 2.5 & 7 \\
1.0 & 0.1 & -0.5 & 5.0 & 1.0 & 5.0 & 2.778 & 11 \\
0.5 & 0.1 & -1.0 & 5.0 & 3.0 & 1.0 & 1.125 & 9 \\
0.5 & 0.1 & -2.0 & 5.0 & 3.0 & 0.5 & 1.125 & 8.5 \\
3.0 & 0.5 & -1.0 & 5.0 & 1.0 & 5.0 & 6.9 & 11 \\
5.0 & 1.0 & -2.0 & 5.0 & 1.0 & 5.0 & 12.5 & 11 \\
\hline
\end{tabular}
\caption{Examples of parameter values giving stable periodic solutions with the coefficient's period $T$.}
\label{table_single_T}
\end{table*}

%%
%% Type II periodic solutions
%%
\subsection{Periodic solutions of type II}\label{PS-II}

\begin{figure*}
    \centering
    \includegraphics[width=0.9\textwidth]{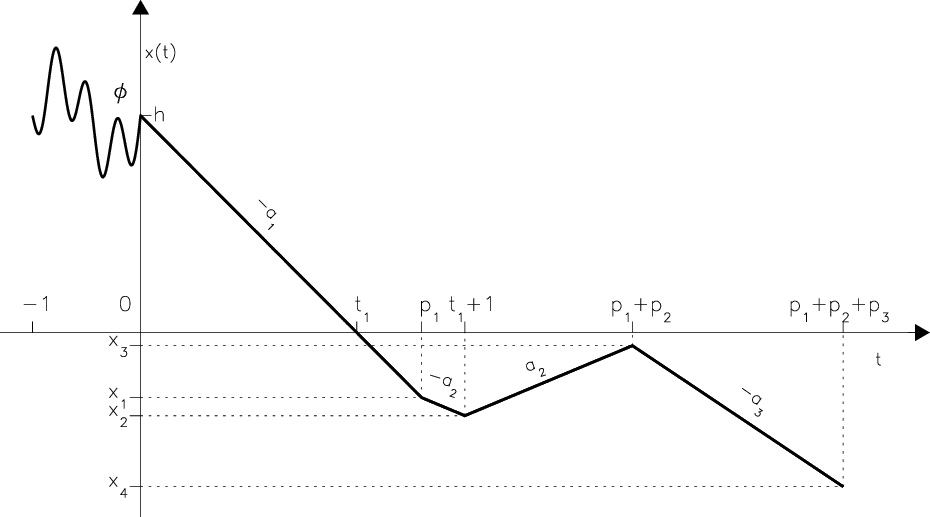}
    \caption{Slowly oscillating piecewise affine solution (Type II) on initial periodic interval $[0,T]$.}
    \label{fig2}
\end{figure*}

Following our previous results in \cite{IvaShel2024_MCA} and based on recent numerical simulations of our particular DDE (\ref{eq1}) we are 
looking next for periodic solutions of the form as shown in Figure~\ref{fig2}. I.e., given an initial function $\varphi\in C$ with 
$\varphi(s)>0\; \forall s\in[-1,0], \varphi(0)=h>0,$ there is only one sign change of the corresponding solution $x=x_{\varphi}(t)$ 
on the first periodic segment $(0,T)$. That is, there is a first zero at $t=t_1$ with $t_1<p_1$ and 
$x(t)>0 \; \forall t\in[0,t_1), x(t)<0\; \forall t\in (t_1,T]$.

%%
%% Fig.2
%%

Define the following values of the corresponding solution $x$ on the interval $(0,T)$: 
$x_1=x(p_1), x_2=x(t_1+1), x_3=x(p_1+p_2), x_4=x(p_1+p_2+p_4).$ 
In order to obtain the solution of the form shown in Figure~\ref{fig2} one must assume that the condition $x_3<0$ is satisfied. 
Due to the piecewise affine form of the solution on the interval $(0,T)$ it is straightforward to calculate the respective values. 
One has the following outcomes:
$$
t_1=\frac{h}{a_1},\; x_1=x(p_1)= h-a_1p_1<0, \;
$$
$$
x_2=x(t_1+1)=\left(1-\frac{a_2}{a_1}\right)h-a_2+p_1(a_2-a_1)
$$
$$
x_3=x(p_1+p_2)=\left(1-2\frac{a_2}{a_1}\right)h-a_1p_1-a_2(2-2p_1-p_2)
$$
\begin{align*}
    x_4=& x\left(p_1+p_2+p_3\right)=x_3-a_3p_3\\
       =& (1-2\frac{a_2}{a_1})h-a_1p_1-a_2(2-2p_1-p_2)-a_3p_3\\
       :=& F_1(h)=mh-b,
\end{align*}
where
\begin{equation}\label{mb}
m=1-2\frac{a_2}{a_1}\;\text{and }\; b=a_1p_1+a_2(2-2p_1-p_2)+a_3p_3>0.
\end{equation}
Likewise, for an initial function $\psi\in C$ satisfying $\psi(s)<0\; \forall s\in[-1,0]$ with $h<0$ one constructs a symmetric 
solution to that shown in Figure~\ref{fig2} and such that $x(t)<0\; \forall t\in[0,t_1)$ and $x(t)>0\; \forall t\in(t_1,T]$. A similar calculation
to the one above yields a map $x_4=F_2(h)=mh+b, h<0,$ where the values $m,b$ are given my the same formulas (\ref{mb}).

It is clear that the fixed points $h_*$ of the composite map $F=F_2\circ F_1$ result in slowly oscillating periodic solutions to DDE 
(\ref{DDE-per}) with the period $2T$. Moreover, such a periodic solution is asymptotically stable if $|m|<1$ and unstable when $|m|>1$. 
The condition $|m|<1$ is equivalent to $a_1>a_2$. The condition $b>0$ is implied from the assumption that $x_3<0$ (first consideration above), 
which is equivalent to $a_1p_1+a_2(2-2p_1-p_2)>0$. The unique cycle of period two of the map $F$ is easily found as $h_*,-h_*$, where
\begin{equation}\label{h2}
h_*=\frac{b}{m+1}=\frac{a_1[a_1p_1+a_2(2-2p_1-p_2)+a_3p_3]}{2(a_1-a_2)}.
\end{equation}

Therefore, we arrive at the following result:
%%
%% Theorem 4.2
%%
\begin{theorem}\label{thm2}
  Suppose that the parameters $a_1,a_2,a_3,p_1,p_2,p_3$ are such that $a_1>a_2$ and the inequality $a_1p_1+a_2(2-2p_1-p_2>0)$ is satisfied. 
Then DDE (\ref{DDE-per}) has an asymptotically stable slowly oscillating periodic solution of period $2T=2(p_1+p_2+p_3)$. The periodic 
solution is generated by the initial function $\phi(s)\equiv h_*, s\in[-1,0],$ where $h_*$ is given by (\ref{h2}).
\end{theorem}

As before, numerical tests confirm existence of these solutions. Examples of parameters, for which Equation~\ref{DDE-per} gives stable slowly-oscillating $2T$-periodic solutions, are provided in Table~\ref{table_double_T}.

\begin{table*}[h!]
\centering
\begin{tabular}{|c|c|c|c|c|c|c|c|c|c|}
\hline
$a_1$ & $a_2$ & $a_3$ & $p_1$ & $p_2$ & $p_3$ & $h^*$ & $T$ & $2T$ \\
\hline
0.5 & 0.1 & -0.1 & 0.5 & 1.0 & 0.5 & 0.1875 & 2 & 4 \\
0.5 & 0.1 & -0.1 & 1.0 & 3.0 & 0.5 & 0.1562 & 4.5 & 9 \\
0.5 & 0.1 & -1.0 & 3.0 & 1.0 & 1.0 & 1.25 & 5 & 10 \\
1.0 & 0.1 & -0.1 & 1.0 & 3.0 & 0.5 & 0.4167 & 4.5 & 9 \\
1.0 & 0.5 & -2.0 & 3.0 & 1.0 & 1.0 & 2.5 & 5 & 10 \\
1.0 & 0.1 & -0.1 & 0.5 & 5.0 & 0.5 & 0.0833 & 6 & 12 \\
1.0 & 0.1 & -0.5 & 0.5 & 3.0 & 0.5 & 0.3056 & 4 & 8 \\
3.0 & 0.1 & -0.1 & 0.5 & 3.0 & 1.0 & 0.7241 & 4.5 & 9 \\
3.0 & 0.1 & -1.0 & 1.0 & 1.0 & 0.5 & 1.759 & 2.5 & 5 \\
3.0 & 0.1 & -2.0 & 0.5 & 5.0 & 0.5 & 1.086 & 6 & 12 \\
5.0 & 0.1 & -0.1 & 0.5 & 1.0 & 0.5 & 1.301 & 2.0 & 4 \\
5.0 & 0.1 & -0.1 & 0.5 & 3.0 & 0.5 & 0.3438 & 4.0 & 8 \\
5.0 & 1.0 & -0.5 & 0.5 & 1.0 & 1.0 & 1.875 & 2.5 & 5 \\
5.0 & 1.0 & -2.0 & 1.0 & 1.0 & 1.0 & 3.75 & 3 & 6 \\
5.0 & 1.0 & -2.0 & 1.0 & 3.0 & 1.0 & 2.5 & 5 & 10 \\
5.0 & 3.0 & -2.0 & 1.0 & 1.0 & 0.5 & 3.75 & 2.5 & 5 \\
10.0 & 0.1 & -0.5 & 0.5 & 3.0 & 1.0 & 2.677 & 4.5 & 9 \\
10.0 & 0.5 & -2.0 & 0.5 & 1.0 & 0.5 & 3.158 & 2 & 4 \\
10.0 & 3.0 & -1.0 & 1.0 & 3.0 & 1.0 & 1.429 & 5 & 10 \\
10.0 & 5.0 & -2.0 & 1.0 & 1.0 & 1.0 & 7.0 & 3 & 6 \\
\hline
\end{tabular}
\caption{Samples of parameter values resulting in stable periodic solutions with the double period $2T$ of the coefficient $a$.}
\label{table_double_T}
\end{table*}
%%
%% Subsection (d)
%%
\subsection{Smoothed continuous nonlinearities}\label{smooth}

Though the piecewise constant DDEs provide a precise and accurate insight into specific features of their dynamics, by exact reduction to interval maps,  
a formal theoretical issue behind their use has always been the fact  that the functions involved are not continuous. 
This leads to the solutions of DDEs being not differentiable at a discrete set of isolated points (as it is can be seen in Figures~\ref{fig1} and \ref{fig2}).
This point can usually be overcome by the procedure of smoothing of the nonlinearities in a small neighbourhood of their discontinuity points in such a way
that the principal dynamics of interest remain intact. This approach of smoothing was successfully used by most of the 
related preceding papers, see, e.g., a sample selection of works \cite{HalLin86,Wal81b}. Regarding DDE (\ref{DDE-per}) the smoothing procedures were 
successfully implemented in our recent publications \cite{IvaShel2024_MCA,IvaLWShel2024}. The simplest approach of making functions $f_0(x)$ and $A_0(t)$
of DDE (\ref{eq1}) continuous is using matching connecting line  segments in a small neighbourhood of their discontinuity points $x=0$ and 
$t=0, p_1, p_2, p_3$. Namely, given small $\delta>0$, define functions $f=f_{\delta}(x)$ and $a=A_{\delta}(t)$ as follows 
\begin{equation}\label{f_d}
    f_{\delta}(x)=
\begin{cases}
+1\; \text{if}\; x\le -\delta\\
-1\; \text{if}\; x\ge \delta\\
-({1}/{\delta}) x\; \text{if}\; x\in(-\delta, \delta),
\end{cases}
\end{equation}
and
\begin{equation}\label{A0d}
A_{\delta}(t)=\begin{cases}
a_3+\frac{a_1-a_3}{2\delta}(t+\delta)\;\text{if}\; t\in(-\delta,\delta)\\
a_1\;\text{if}\; t\in[\delta,p_1-\delta]\\
a_1+\frac{a_2-a_1}{2\delta}[t-(p_1-\delta)]\;\text{if}\; t\in(p_1-\delta,p_1+\delta)\\
a_2\;\text{if}\; t\in[p_1+\delta,p_1+p_2-\delta]\\
a_2+\frac{a_3-a_2}{2\delta}[t-(p_2-\delta)]\;\\
\text{if}\; t\in(p_1+p_2-\delta,p_1+p_2+\delta)\\
a_3\;\text{if}\; t\in[p_1+p_2+\delta,p_1+p_2+p_3-\delta]\\
a_3+\frac{a_1-a_3}{2\delta}[t-(p_1+p_2+p_3-\delta)]\;\\
\text{if}\; t\in(p_1+p_2+p_3-\delta,p_1+p_2+p_3+\delta)\\
\text{periodic extension on}\; \mathbb{R}\; \text{outside interval}\;\\ 
[0,T), T=p_1+p_2+p_3.
\end{cases}
\end{equation}

%\textcolor{red}{ ... to describe the fact that the periodic dynamics persist for small $\delta\ge0,$ use the results from our MCA 2024 paper (as a theorem).}

The smoothing of discontinuities procedure done by using continuous functions (\ref{f_d}) and (\ref{A0d}) has an advantage compared with other options. 
It turns out that in this case the new smoothed solutions (of $C^1$-class) of delay differential equation (\ref{P-C dde}) maintain exactly the same shape and 
position outside of a small (of order $\delta$) neighbourhood of the previous "corner points". Those are the $t$-values $p_1, t_1+1, p_1+p_2, p_1+p_2+p_3$ for the
solution shown in Figure~\ref{fig1} (similar $t$-values are shown in Figure~\ref{fig2}). Due to the straight line connecting segments, the corresponding smooth solutions of DDEs are 
parabolas around the previously ``corner points'' that are symmetric with respect to those values. 
This implies the fact that the maps $F, F_1, F_2$  defined in subsections 4(b) and 4(c) do not change under such small $\delta$-perturbation. 
This leads to the same fixed points $h_*$ and their types of stability. This result of the ``no change in the map'' was proved in \cite{IvaShel2024_MCA}. 
That proof does not change in main steps and straightforwardly extends to the case of functions $f_\delta$ and $A_\delta$ in DDE  (\ref{P-C dde}). See relevant details 
in our paper \cite{IvaShel2024_MCA}: page 11 for the comparison graphs of  ``corner'' and ``smoothed'' solutions, and Theorem 4 on page 12 for the equivalence of 
the dynamics. 

We should note that, in general, the type of the monotone smoothing connection on the $\delta$-neighbourhood of the discontinuity set of $f$ and $a$ does not really 
matter. The one-dimensional dynamics, its fixed points, and their stability will persist and depend continuously on $\delta\ge0$ as $\delta\to0$ (with $\delta=0$ 
being the piecewise constant discontinuous functions $f_0(x)$ and $A_0(t)$). This general result is proved in our recent paper \cite{IvaLWShel2024} for equation (\ref{DDE-per}). See Theorems 3.5 and 3.6 there for details. 

\section{Conclusion and Discussion}\label{Concl}

This paper demonstrates, by presenting explicit examples with precise calculations, that delay differential equations of the form (\ref{DDE-per}) 
with $T$-periodic coefficients and parameters can possess nontrivial periodic solutions with the same period $T$ or double of it. With respect to the applied 
problems and the natural phenomena behind them, that are modelled by  the equations,  this outcome  can be interpreted that the periodic input on equations 
leads to the existence of periodic motions with the same period or a multiple of it. 
Therefore, the general problem that one can pose for the DDE (\ref{DDE-per}) is to derive the sufficient conditions for the 
existence of nontrivial $T$-periodic solutions under the $T$-periodicity assumptions on the coefficients $\mu(t), a(t)$ and the delay $\tau(t)$. 
The first natural step for solving this problem can be, based on the properties of the prototype autonomous models, to derive such conditions when 
the decay rate $\mu(t)=\mu_0>0$ and the multiplier periodic coefficients are positive, and the delay $\tau(t)=\tau_0>$ is a constant.
This can be relatively easy done, in our opinion, for model examples of the type we treated in this paper, and following the related calculations already in place  
(as they are done in our recent work \cite{IvaLWShel2024} with the constant decay term $\mu>0$). However, the related exact calculations are much more involved 
than those in the present paper; such a new research work can be a subject of a separate paper.

%\textcolor{blue}{It makes sense to talk here about the general problem of existence of periodic solutions to DDE (\ref{DDE-per}), e.g. in the sense of the following sample aspects
%\begin{itemize}
%    \item 
%    Deriving sufficient conditions for the existence of such periodic solutions under general assumption of a (negative/positive) feedback assumption on nonlinearity $f$ and the periodicity of the decay and multiplicative coefficients $\mu(t)$ and $a(t)$.
%    \item 
%    Discuss possible transition in dynamics from the autonomous DDEs to those with periodic coefficients. Perhaps limit the discussion to the simplest case when $\mu=0$, $f$ is of negative feedback, and $a(t)>0$ is parameter dependent with the amplitude of the periodic oscillations increasing from zero (autonomous case) to large ones.
%    \item 
%    It would be nice to have some numerical insight (even limited preliminary calculations) to support some of the features discussed in the two items above.
%\end{itemize}}

%\cite{Had79} \cite{MacGla77}

\section{Acknowledgement}

This work was initiated during A. Ivanov's visit to the University of Sydney as a research fellow at the Sydney Mathematical Research 
Institute within their 
International Visitor Program (November 2024 - January 2025). The research program also included collaborative activities at Deakin University 
(G. Beliakov) and Flinders University (S. Shelyag). The research was completed during a follow-up visit to Flinders University in August 2025. A. Ivanov expresses his gratitude for the support and accommodation at all three universities.

%\disclaimer{Insert disclaimer text here.}

%%%%%%%%%% Insert bibliography here %%%%%%%%%%%%%%

\bibliographystyle{alpha}
\bibliography{refs}

\begin{thebibliography}{adHW83}

\bibitem[adHW83]{AdHW}
U.~an~der Heiden and H.-O. Walther.
\newblock Existence of chaos in control systems with delayed feedback.
\newblock {\em J. Differential Equations}, 47:273--293, 1983.

\bibitem[dHM82]{AdHM}
U.~An der Heiden and M.C. Mackey.
\newblock The dynamics of production and destruction: Analytic insight into complex behavior.
\newblock {\em J. Math. Biol.}, 16:75--101, 1982.

\bibitem[ea95]{DieSvGSVLWal95}
O.~Diekmann et~al.
\newblock Delay equations: Complex, functional, and nonlinear analysis.
\newblock {\em Springer-Verlag, Ser.: Applied Mathematical Sciences}, 110, 1995.

\bibitem[Ern09]{Ern09}
T.~Erneux.
\newblock {\em Applied Delay Differential Equations}, volume~3 of {\em Ser.: Surveys and Tutorials in the Applied Mathematics Sciences}.
\newblock Springer Verlag, 2009.

\bibitem[GBN80]{GurBlyNis80}
W.S.C. Gurney, S.P. Blythe, and R.M. Nisbet.
\newblock Nicholson's blowflies revisited.
\newblock {\em Nature}, 287:17--21, 1980.

\bibitem[GM88]{GlaMac88}
L.~Glass and M.C. Mackey.
\newblock {\em From Clocks to Chaos: Rhythms of Life}.
\newblock Princeton University Press, 1988.

\bibitem[Had79]{Had79}
K.P. Hadeler.
\newblock Delay equations in biology.
\newblock {\em Springer-Verlag, Ser.: Applied Mathematical Sciences}, 730:139--156, 1979.

\bibitem[HL86]{HalLin86}
J.K. Hale and X.B. Lin.
\newblock Examples of transverse homoclinic orbits in delay equations.
\newblock {\em Nonlinear Anal.}, 10:693--709, 1986.

\bibitem[HL93]{HalSVL93}
J.K. Hale and S.M.~Verduyn Lunel.
\newblock Introduction to functional differential equations.
\newblock {\em Springer Applied Mathematical Sciences}, 99, 1993.

\bibitem[ILWS24]{IvaLWShel2024}
Anatoli Ivanov, Bernhard Lani-Wayda, and Sergiy Shelyag.
\newblock Periodic solutions of a delay differential equation with a periodic multiplier, 2024.

\bibitem[IS92]{IvaSha92}
A.F. Ivanov and A.N. Sharkovsky.
\newblock Oscillations in singularl perturbed delay equations.
\newblock {\em Dynamics Reported (New Series)}, 1:165--224, 1992.

\bibitem[IS23]{IvaShel2023}
Anatoli Ivanov and Sergiy Shelyag.
\newblock Stable periodic solutions in scalar periodic differential delay equations.
\newblock {\em Archivum Mathematicum}, (1):69–76, 2023.

\bibitem[IS24a]{IvaShel2024}
Anatoli Ivanov and Sergiy Shelyag.
\newblock Explicit periodic solutions in a delay differential equation, 2024.

\bibitem[IS24b]{IvaShel2024_MCA}
Anatoli Ivanov and Sergiy Shelyag.
\newblock Periodic solutions in a simple delay differential equation.
\newblock {\em Mathematical and Computational Applications}, 29(3), 2024.

\bibitem[IS25]{IvaShel2025}
Anatoli Ivanov and Sergiy Shelyag.
\newblock {\em Delay Differential Equations with Periodic Coefficients: A Numerical Insight}, pages 433--439.
\newblock Springer Nature Switzerland, Cham, 2025.

\bibitem[JFPC78]{PerMalCou78}
C.~P.~Malta J.~F.~Perez and F.~A.~B. Coutinho.
\newblock Qualitative analysis of oscillations in isolated populations of flies.
\newblock {\em J. Theoret. Biol.}, 71:505--514, 1978.

\bibitem[KM99]{KolMys99}
V.~Kolmanovskii and A.~Myshkis.
\newblock {\em Theory and Applications of Functional Differential Equations}.
\newblock Ser.: Mathematics and Its Applications. Kluwer Academic Publishers, 1999.

\bibitem[Kua93]{Kua93}
Yang Kuang.
\newblock Delay differential equations with applications in population dynamics.
\newblock {\em Academic Press Inc., Series: Mathematics in Science and Engineering}, 191:398, 1993.

\bibitem[MG77]{MacGla77}
M.C. Mackey and L.~Glass.
\newblock Oscillation and chaos in physiological control systems.
\newblock {\em Science}, 197:287--289, 1977.

\bibitem[Pet80]{Pet80}
H.~Peters.
\newblock Comportment chaotique d’une equation retardee.
\newblock {\em C. R. Acad. Sci. Paris Ser. I. Math.}, 290:1119--1122, 1980.

\bibitem[Smi11]{Smi11}
H.~Smith.
\newblock {\em An Introduction to Delay Differential Equations with Applications to the Life Sciences}, volume~57.
\newblock Springer-Verlag, Series: Texts in Applied Mathematics, 2011.

\bibitem[Wal81]{Wal81b}
H.-O. Walther.
\newblock Homoclinic solution and chaos in $\dot x(t)=f(x(t-1))$.
\newblock {\em Nonlinear Analysis, TMA}, 5:775--788, 1981.

\end{thebibliography}

\end{document}